\documentclass[11pt,a4paper]{article}

\usepackage{amsmath}
\usepackage{amsthm}
\usepackage{amsfonts}

\parskip\smallskipamount
\newtheorem{theorem}{Theorem}
\newtheorem{lem}{Lemma}
\newtheorem{prop}{Proposition}
\theoremstyle{remark}
\newtheorem{rem}{Remark}

\renewcommand{\Pr}{\mathbf{P}}
\newcommand{\E}{\mathbf{E}}
\newcommand{\R}{\mathbb{R}}
\newcommand{\Rp}{\mathbb{R}_+}
\newcommand{\I}{\textbf{I}}
\newcommand{\ov}[1]{\overline{#1}}
\newcommand{\LL}[1]{\mathcal{L}_#1}
\renewcommand{\SS}[1]{\mathcal{S}_#1}
\renewcommand{\phi}{\varphi}
\renewcommand{\epsilon}{\varepsilon}
\newcommand{\phg}[1]{\hat{\phi}_{#1}}

\parskip \smallskipamount

\title{On the exact distributional asymptotics for the supremum
  of a random walk with increments in a class of light-tailed
  distributions}
\author{Stan Zachary and Serguei Foss\\
  \emph{Heriot-Watt University, Edinburgh}\\
  \emph{and} \\
  \emph{Institute of Mathematics,  Novosibirsk}}
\date{\today}

\begin{document}

\maketitle






\vspace{-1ex}

\begin{quotation}\small
  We study the distribution of the maximum~$M$ of a random walk whose
  increments have a distribution with negative mean and belonging, for
  some $\gamma>0$, to a subclass of the class~$\SS{\gamma}$---see, for
  example, Chover, Ney, and Wainger (1973).  For this subclass we give
  a probabilistic derivation of the asymptotic tail distribution of
  $M$, and show that extreme values of $M$ are in general attained
  through some single large increment in the random walk near the
  beginning of its trajectory.  We also give some results concerning
  the ``spatially local'' asymptotics of the distribution of $M$, the
  maximum of the stopped random walk for various stopping times, and
  various bounds.
\end{quotation}

\vspace{-1ex}

\section{Introduction}
\label{sec:introduction}

For any distribution function~$F$ on $\R$ define the
function~$\phi_F\colon\Rp\to\Rp\cup\{\infty\}$ by
\begin{displaymath}
  \phi_F(\alpha) := \int_{-\infty}^\infty e^{\alpha x}dF(x) = \E e^{\alpha X}
\end{displaymath}
where the random variable~$X$ has distribution function~$F$.  Denote
also by $\ov{F}$ the tail distribution function given by
$\overline{F}(x)=1-F(x)$.  Let $F^{*n}$ denote the
$n$-fold convolution of $F$ with itself.

As usual, we shall say that a distribution function~$F$ on
$\R$ belongs to the class~$\LL{\gamma}$, $\gamma\ge0$, if and only if
\begin{equation}
  \label{eq:1}
    \text{$\ov{F}(x)>0$ for all $x$},\qquad
    \lim_{x\to\infty}\frac{\ov{F}(x-k)}{\ov{F}(x)} = e^{\gamma k},
    \quad\text{for all fixed $k$}.
\end{equation}
Note that the class $\LL{0}$ is the usual class of \emph{long-tailed}
distributions.  As in that case, for $F\in\LL{\gamma}$,
the convergence above is necessarily uniform for all $k$ in
any compact interval.  Further, if $F\in\LL{\gamma}$, then $\gamma$ is
uniquely defined by
$\gamma=\sup\{\alpha\colon{}\phi_F(\alpha)<\infty\}$.  In this case we
shall find it convenient to write~$\phg{F}$ for $\phi_{F}(\gamma)$,
which may be either finite or infinite.  Note also that if
$F\in\LL{\gamma}$ and $\phg{F}<\infty$, then
\begin{equation}
  \label{eq:2}
  \ov{F}(x) = o(e^{-\gamma x}),
  \qquad\text{as $x\to\infty$}.
\end{equation}

A distribution function~$F$ on $\R$ belongs to the
class~$\SS{\gamma}$, $\gamma\ge0$, if and only if $F\in\LL{\gamma}$,
$\phg{F}<\infty$, and
\begin{equation}
  \label{eq:3}
  \ov{F^{*2}}(x) \sim 2\phg{F}\ov{F}(x),
  \qquad\text{as $x\to\infty$}
\end{equation}
(where, for any two positive functions~$f_1$, $f_2$ on $\R$, we write
$f_1(x)\sim{}f_2(x)$ as $x\to\infty$ if $f_1(x)/f_2(x)\to1$ as
$x\to\infty$).  We remark that $\SS{\gamma}$ is sometimes defined
first for distribution functions~$F$ on $\Rp$, and the class then
extended to $F$ on the whole real line
by requiring $F\I_{\Rp}\in\SS{\gamma}$, where
$\I_{\Rp}$ is the indicator function of $\Rp$.  However, this is
unnecessary, essentially because of the condition~$F\in\LL{\gamma}$.

Let $M$ be the maximum of a random walk whose increments have
distribution function~$F$ with negative mean.  For $F$ belonging to
the well-known class~$\SS{0}$ of subexponential distributions, the
asymptotic distribution of the tail of $M$ was given by Pakes (1975)
and Veraverbeke (1977), see also the earlier result by Borovkov
(1976) for the subclass of distributions with regularly varying tails.
In this case extreme values of $M$ are in
general attained through some single large increment in the random
walk---see Zachary (2004).  This is the ``principle of a single big
jump''.  Further, conditional on $M$ exceeding $x$, the probability
that it does so within any given finite time tends to zero as $x$
tends to infinity.

For $F\in\SS{\gamma}$ for some $\gamma>0$, it necessary to distinguish
various subclasses of $\SS{\gamma}$.  In the case where $\phg{F}>1$,
or where $\phg{F}=1$ and $\phi_{F}'(\gamma)<\infty$, the asymptotic
distribution of $M$ is as given by the classical Cram\'{e}r-Lundberg
theory; in particular trajectories leading to extreme values of $M$
are typically approximately linear prior to the time at which $M$ is
attained.  However, our interest in the present paper is primarily in
that subclass of $\SS{\gamma}$ for which $\phg{F}<1$.  Here the
asymptotic distribution of the tail of $M$ was given by Borovkov
(1976, Chapter~4) for the case where the function $G$ given by
$G(x)=e^{\gamma{}x}\ov{F}(x)$ is regularly varying.  Further, for $G$
regularly varying, a considerably more general investigation of sample
path behaviour associated with extreme values of $M$ was given by
Borovkov and Borovkov (2004).  The result of Borovkov was extended
from the regularly varying to the general case by Veraverbeke (1977).
Bertoin and Doney (1996) showed that there is a gap in Veraverbeke's
proof and proposed their own corrected version.  All these proofs are
analytic.  Our present aim is to derive in this case the asymptotic
form of $\Pr(M>x)$, as $x\to\infty$, using arguments which are
entirely probabilistic.  We also give a number of corollaries and
associated results, which in particular provide further insight into the
behaviour of the sample paths leading to extreme values of $M$.  These
show that for the case $\gamma>0$, $\phg{F}<1$, the principle of a
single big jump again holds, i.e.\ extreme values of $M$ are in
general again attained through some single large increment in the
random walk.  However, in this case such extreme values are likely to
occur close to the start of this process---see Remark~\ref{rem:bj2}
below for a precise statement, and again Borovkov and Borovkov (2004)
for more on the case where $G$ is regularly varying.  We also present
results for the (spatially) local asymptotics for $M$, the asymptotics
for the maximum of the stopped random walk (for various stopping
times), and we give some bounds.

Finally in this section we give some general properties of the 
class~$\SS{\gamma}$, $\gamma\ge0$, which will be required
subsequently.
For $F\in\LL{\gamma}$, $\gamma\ge0$, with
$\phg{F}<\infty$, it is straightforward to show that~\eqref{eq:3} is
equivalent to the requirement that
\begin{equation}
  \label{eq:4}
  \lim_{x\to\infty}\frac{1}{\ov{F}(x)}\int_{h(x)}^{x-h(x)}
  dF(y)\ov{F}(x-y) = 0,
\end{equation}
for any positive function $h$ such that
\begin{equation}
  \label{eq:5}
  h(x)\le{}x/2 \text{ for all $x$},
  \qquad \lim_{x\to\infty}h(x)=\infty.
\end{equation}
The class $\SS{\gamma}$ is thus indeed a natural generalisation of the
class~$\SS{0}$ of subexponential distributions (recall also that
$\phg{F}=1$ for $F\in\SS{0}$).  As in the case~$\gamma=0$, it is not
easy to give an example of a distribution~$F$ which belongs to
$\LL{\gamma}$ with $\phg{F}<\infty$, but not to $\SS{\gamma}$, and all
such examples are more-or-less artificial.


It is clear that, for $\gamma\ge0$, the class $\SS{\gamma}$ is closed
under tail-equivalence (i.e., if $F\in\SS{\gamma}$ and
$\ov{F}(x)\sim{}c\ov{G}(x)$ with $0<c<\infty$,
then $G\in\SS{\gamma}$).  More generally, we have the result given by
Proposition~\ref{prop:1} below, the proof of which follows from
Lemma~5.1 of Pakes (2004) analogously to the proof of Lemma~5.2 of
that paper.

\begin{prop}\label{prop:1}
  Suppose that $F\in\SS{\gamma}$ for some
  $\gamma\ge0$.  Suppose further that, for $i=1,\dots,n$, the
  distribution function~$F_i$ is such that
  $\ov{F_i}(x)\sim{}c_i\ov{F}(x)$ as $x\to\infty$ for some $c_i\ge0$
  (where in the case $c_i=0$ we understand this to mean that
  $\ov{F_i}(x)=o(\ov{F}(x))$ as $x\to\infty$).  Then
  $\phg{F_i}<\infty$ for all $i=1,\dots,n$ and the convolution
  $F_1*\dots*F_n$ satisfies
  \begin{displaymath}
    \ov{F_1*\dots*F_n}(x)
    \sim
    \prod_{i=1}^n\phg{F_i}\sum_{i=1}^n\frac{c_i}{\phg{F_i}}
    \;\ov{F}(x),
    \qquad\text{as $x\to\infty$,}
  \end{displaymath}
  Further, if $\sum_{i=1}^n{}c_i>0$ then $F_1*\dots*F_n\in\SS{\gamma}$.
\end{prop}
In particular we have the generalisation of the property~\eqref{eq:3} above,
that if $F\in\SS{\gamma}$, then, for all $n\ge1$,
\begin{displaymath}
  \ov{F^{*n}}(x) \sim n\phg{F}^{n-1}\ov{F}(x),
  \qquad\text{as $x\to\infty$}.
\end{displaymath}

For future reference, we note that, for $\gamma\ge0$, the
class~$\SS{\gamma}$ may be extended to include distributions~$F$ with
support on
$\R\cup\{-\infty\}$, i.e.\ we may allow the possibility of strictly
positive mass at~$-\infty$.  For $\gamma>0$, Proposition~\ref{prop:1}
remains unchanged.  This is also true for $\gamma=0$, provided that
here, for any $F\in\SS{0}$, we take $\phg{F}=\ov{F}(-\infty)$.

\section{Random walks with negative drift}
\label{sec:random-walks-with}

Let $\{\xi_n\}_{n\ge1}$ be independent identically distributed random
variables with distribution function~$F$ on $\R$.  Let $S_0=0$,
$S_n=\sum_{i=1}^n\xi_i$ for $n\ge1$.  Let
$M_n=\max_{0\le{}i\le{}n}S_i$ for $n\ge0$ and let $M=\sup_{n\ge0}S_n$.

Suppose that $F\in\LL{\gamma}$ for some $\gamma>0$ and that
$\phg{F}<1$.  Since $\phi_F(0)=1$ and $\phi_F$ is convex on
$[0,\gamma]$, it follows that $\phi_F'(0)<0$ and so $F$ has a negative
mean, which, we note, may be $-\infty$.  It thus follows that
$\Pr(M<\infty)=1$.  We are interested in the asymptotic distribution
of $\Pr(M>x)$ as $x\to\infty$.

The following lemma gives a useful preliminary result.
\begin{lem}\label{lem:1}
  Suppose that $F\in\LL{\gamma}$ for some $\gamma>0$ and that
  $\phg{F}<1$.  Then
  \begin{equation}
    \label{eq:33}
    1 \le \E e^{\gamma M} \le \frac{1}{1-\phg{F}}
  \end{equation}
  and therefore
  \begin{displaymath}
  \Pr(M > x) = o(e^{-\gamma x}), \qquad\text{as $x\to\infty$.} 
  \end{displaymath}
\end{lem}

\begin{proof}
  Observe that
  \begin{displaymath}
    1 \le e^{\gamma M} \le \sum_{n=0}^{\infty} e^{\gamma S_n}.
  \end{displaymath}
  The result~\eqref{eq:33} now follows by taking expectations and
  noting that $\E e^{\gamma S_n}=\phg{F}^n$.
\end{proof}

For any stopping time $\sigma$ define, on the set $\sigma<\infty$,
\begin{displaymath}
  M^\sigma = \sup_{m\ge0} S_{\sigma+m} - S_\sigma.
\end{displaymath}
(In particular, for any finite $n$, $M^n=\sup_{m\ge0}S_{n+m}-S_n$.)
For each $x\ge{}a\ge0$, and each $n\ge1$, define the event
\begin{displaymath}
  A_n^{a,x} = \{M_{n-1} \le a,\, S_n > x\}.
\end{displaymath}
Note that for any $a$ and $x$ as above, the events $A_n^{a,x}$,
$n\ge1$, are disjoint.  The following lemma gives a lower bound.

\begin{lem}
  \label{lem:lb}
  Suppose that $F\in\LL{\gamma}$ for some $\gamma>0$ and
    that $\phg{F}<1$.  Then, given $\epsilon>0$, there exists $a>0$
    such that
    \begin{displaymath}
      \liminf_{x\to\infty}\frac{\Pr(M>x)}{\ov{F}(x)}
      \ge
      \liminf_{x\to\infty}\frac{1}{\ov{F}(x)}
      \Pr\biggl(M>x,\,\bigcup_{n\ge1}A_n^{a,x-a}\biggr)
      \ge
      (1-\epsilon)\frac{\E e^{\gamma M}}{1-\phg{F}}.
    \end{displaymath}
\end{lem}

\begin{proof}
  Observe first that, from the definition~\eqref{eq:1} and the
  monotonicity of $\ov{F}$, it follows straightforwardly that, given any
  $\gamma'\in(0,\gamma)$, there exists $x_0$ such that, for all $x\ge
  x_0$ and all $t\le-1$,
  \begin{align}
    \ov{F}(x-t) & \le \ov{F}(x)e^{\gamma t}e^{(\gamma'-\gamma) t}
    \nonumber\\
    & = \ov{F}(x)e^{\gamma' t}. \label{eq:10}
  \end{align}
  We thus have, for any $a\in\R$ and any $n\ge1$,
  \begin{align}
    \Pr(S_{n-1}\le a,\,S_n>x)
    & = \int_{-\infty}^a \Pr(S_{n-1}\in dt)\ov{F}(x-t)\label{eq:11}\\
    & = (1+o(1))\ov{F}(x) \int_{-\infty}^a \Pr(S_{n-1}\in dt)
    e^{\gamma t}\label{eq:12}
  \end{align}
  as $x\to\infty$.  This follows from \eqref{eq:1} and the dominated
  convergence theorem, where \eqref{eq:10} and \eqref{eq:1}, with the
  asserted uniformity, are used to bound the integrand respectively in
  the regions $(-\infty,-1)$ and $[-1,a]$.


  Now fix $a>0$.  For each $n\ge1$, and for $x\ge a$,
  \begin{align}
    \Pr(A_n^{a,x})
    & = \Pr(S_{n-1}\le a,\,S_n>x)
    - \Pr(M_{n-2}>a,\, S_{n-1}\le a,\,S_n>x)\nonumber\\
    & \ge \Pr(S_{n-1}\le a,\,S_n>x)
    - \Pr(M_{n-2}>a)\ov{F}(x-a)\nonumber\\
    & \ge \Pr(S_{n-1}\le a,\,S_n>x)
    - \Pr(M>a)\ov{F}(x-a)\nonumber\\
    & = (1+o(1))\lambda(n,a)\ov{F}(x) \label{eq:13}
  \end{align}
  as $x\to\infty$, from~\eqref{eq:1} and \eqref{eq:12}, where
  \begin{equation}
    \label{eq:14}
    \lambda(n,a)
    = \int_{-\infty}^a \Pr(S_{n-1}\in dt)e^{\gamma t} - \Pr(M>a)e^{\lambda a}.
  \end{equation}
  For any $n\ge1$, $x\ge2a$,
  \begin{align}
    \Pr\left(M>x,\,A_n^{a,x-a}\right)
    & \ge \int_0^a \Pr\left(A_n^{a,x-t},\,M^n\in dt\right) \nonumber\\
    & = \int_0^a \Pr\left(A_n^{a,x-t}\right) \Pr(M\in dt) \label{eq:15}\\
    & \ge (1+o(1)) \lambda(n,a) \int_0^a \Pr(M\in dt)\ov{F}(x-t)
    \label{eq:16}\\
    & = (1+o(1)) \lambda(n,a)\ov{F}(x)\int_0^a \Pr(M\in dt)e^{\gamma t}
    \nonumber
  \end{align}
  as $x\to\infty$, where~\eqref{eq:15} follows since, for each $t$,
  the event~$A_n^{a,x-t}$ and the random variable~$M^n$ are
  independent and the latter is equal in distribution to $M$, and
  where \eqref{eq:16} follows from \eqref{eq:13}.  Since also the
  events~$A_n^{a,x-a}$, $n\ge1$, are disjoint, we have
  \begin{equation}
    \label{eq:17}
    \liminf_{x\to\infty}\frac{1}{\ov{F}(x)}
    \Pr\biggl(M>x,\,\bigcup_{n\ge1}A_n^{a,x-a}\biggr)
    \ge
    \sum_{n=1}^N \lambda(n,a)\int_0^a \Pr(M\in dt)e^{\gamma t}
  \end{equation}
  for all $N$.
  Recall that $\phg{F}<1$.    From~\eqref{eq:14}, the
  definition of $\phg{F}$, and Lemma~\ref{lem:1} we see that
  \begin{equation}
    \label{eq:18}
    \lim_{N\to\infty}\lim_{a\to\infty}\sum_{n=1}^N\lambda(n,a)
    = \sum_{n\ge1}\phg{F}^{n-1} = \frac{1}{1-\phg{F}},
  \end{equation}
  while, also
  $\lim_{a\to\infty}\int_0^a\Pr(M\in{}dt)e^{\gamma{}t}=\E{}e^{\gamma{}M}$,
  so that the required result now follows from~\eqref{eq:17}.
\end{proof}

We now give the companion upper bound, which requires the stronger
condition that $F\in\SS{\gamma}$.

\begin{lem}
  \label{lem:ub}
  Suppose that $F\in\SS{\gamma}$ for some $\gamma>0$ and that
  $\phg{F}<1$.  Then
  \begin{displaymath}
    \limsup_{x\to\infty}\frac{\Pr(M>x)}{\ov{F}(x)} \le
    \frac{\E e^{\gamma M}}{1-\phg{F}}.
  \end{displaymath}
\end{lem}

\begin{proof}
  For any sequence of events $\{E_n\}$ we make the
  convention: $\min\{n\ge1:\I(E_n)=1\}=\infty$ if $\I(E_n)=0$ for all
  $n$.  Since $\phg{F}<1$, the distribution~$F$ has a negative mean
  (which, as previously remarked, may be $-\infty$).  Hence there
  exists~$c<0$ such that
  \begin{equation}
    \label{eq:19}
    \limsup_{n\to\infty}S_n/n < c \quad\text{a.s.}
  \end{equation}
  Given also $R>0$, define renewal times
  $0\equiv\tau_0<\tau_1\le\tau_2\le\dots$ for the process~$\{S_n\}$ by
  \begin{displaymath}
    \tau_1 = \min\{n\ge1:S_n>R+nc\} \le \infty,
  \end{displaymath}
  and, for $m\ge2$,
  \begin{align*}
    \tau_m & = \infty, \quad\text{if $\tau_{m-1}=\infty$},\\
    \tau_m & = \tau_{m-1}
    + \min\{n\ge1:S_{\tau_{m-1}+n}-S_{\tau_{m-1}}>R+nc\},
    \quad\text{if $\tau_{m-1}<\infty$}.
  \end{align*}
  Observe that, for each $m\ge1$, conditional on the event
  $\tau_{m-1}<\infty$, the distribution of
  $(\tau_m-\tau_{m-1},\,S_{\tau_m}-S_{\tau_{m-1}})$ is otherwise
  independent of $\mathcal{F}_{\tau_{m-1}}$ (where for each $n\ge1$
  the $\sigma$-algebra~$\mathcal{F}_n$ is that generated by the process
  $\{S_k\}_{k\le{}n}$) and is equal to that of $(\tau_1,\,S_{\tau_1})$.
  In particular $\Pr(\tau_m<\infty)=\delta^m$ where
  \begin{equation}
    \label{eq:20}
    \delta\equiv\Pr(\tau_1<\infty)\to0, \qquad\text{as $R\to\infty$},
  \end{equation}
  from (\ref{eq:19}).  Define also $S_\infty=-\infty$.

  Since the conditions of Lemma~\ref{lem:lb} are also satisfied here,
  it follows from~\eqref{eq:12} that, for any $n\ge1$ and any $a$,
  \begin{equation}
    \label{eq:21}
    \Pr(S_{n-1}\le a,\,S_n>x) \le (1+o(1))\phg{F}^{n-1}\ov{F}(x),
    \qquad\text{as $x\to\infty$}.
  \end{equation}
  Fix also $\gamma'\in(0,\gamma)$ and note that, from \eqref{eq:10}
  and \eqref{eq:11}, for any $n\ge1$ and any $a\le-1$,
  \begin{align}
    \Pr(S_{n-1}\le a,\,S_n>x)
    & \le (1+o(1))\ov{F}(x) \int_{-\infty}^a \Pr(S_{n-1}\in dt) e^{\gamma' t}
    \nonumber\\
    & \le (1+o(1))e^{\gamma' a}\ov{F}(x)
    \qquad\text{as $x\to\infty$},    \label{eq:22}
  \end{align}
  uniformly in $n$ and in $a\le-1$.

  It follows from \eqref{eq:21} and \eqref{eq:22} that, for any $N$
  such that $R+Nc\le-1$, as $x\to\infty$,
  \begin{align}
    \Pr(S_{\tau_1}>x) & = \sum_{n\ge1} \Pr(\tau_1=n,\,S_n>x) \nonumber\\
    & \le \sum_{n\ge1} \Pr(S_{n-1}\le R+(n-1)c,\,S_n>x) \nonumber\\
    & \le (1+o(1))\overline{F}(x)
    \left(
      \sum_{n=1}^N\phg{F}^{n-1}
      + \sum_{n>N}e^{\gamma'(R+(n-1)c)}
    \right). \label{eq:23}
  \end{align}
  Since $c<0$, letting $N\to\infty$, we obtain
  \begin{equation}
    \label{eq:24}
    \Pr(S_{\tau_1}>x) \le
    (1+o(1))\frac{\overline{F}(x)}{1-\phg{F}}
    \qquad\text{as $x\to\infty$}.
  \end{equation}
  
  We now show that, for sufficiently large $R$ and $x$,
  $\sum_{m\ge1}\Pr(S_{\tau_m}>x)$ is comparable to
  $\Pr(S_{\tau_1}>x)$.  Define $\Phi_R=\E e^{\gamma S_{\tau_1}}$.
  Then
  \begin{align}
    \Phi_R & = \sum_{n\ge1} \E\left(\I_{\{\tau_1=n\}} e^{\gamma S_n} \right)
    \nonumber \\
    & \le \sum_{n\ge1} \E\left(\I_{\{S_n>R+nc\}} e^{\gamma S_n} \right)
    \nonumber \\
    & \to 0 \qquad\text{as $R\to\infty$},
    \label{eq:25}
  \end{align}
  by the dominated convergence theorem, since $\phg{F}<1$ and, for
  each $n$,
  \begin{displaymath}
    \E\left(\I_{\{S_n>R+nc\}} e^{\gamma S_n} \right) \le \E e^{\gamma S_n}
    = \phg{F}^n.
  \end{displaymath}
  For each $m\ge1$, the distribution of the random
  variable~$S_{\tau_m}$ is the $m$-fold convolution of the
  distribution of $S_{\tau_1}$ and belongs to $\SS{\gamma}$.  Hence,
  from Proposition~\ref{prop:1} (which, as already noted, extends to
  distributions in $\SS{\gamma}$ with positive mass at $-\infty$),
  \begin{equation}
    \label{eq:26}
    \Pr(S_{\tau_m} > x) = (1+o(1))m\Phi_R^{m-1}\Pr(S_{\tau_1} > x),
    \qquad\text{as $x\to\infty$}
  \end{equation}
  We wish to use again the dominated convergence theorem to
  obtain the corresponding asymptotic result for
  $\sum_{m\ge1}\Pr(S_{\tau_m} > x)$.  To do so we work for the
  moment with conditional distributions.  The distribution of
  $S_{\tau_m}$ conditional on $\tau_m<\infty$ is the $m$-fold
  convolution of the distribution of $S_{\tau_1}$ conditional on
  $\tau_1<\infty$. Further
  \begin{displaymath}
    \E\left(e^{\gamma{}S_{\tau_1}}\,|\,\tau_1<\infty\right)
    = \delta^{-1}\Phi_R
  \end{displaymath}
  (where $\delta$ is as given by \eqref{eq:20}).  It therefore follows
  from Lemma~5.3 of Pakes (2004) that, given $\epsilon>0$, there
  exists a constant~$K$ such that, for all $m\ge1$ and all $x\ge 0$,
  \begin{align*}
    \Pr(S_{\tau_m} > x)
    & = \delta^m \Pr(S_{\tau_m} > x\,|\,\tau_m<\infty)\\
    & \le \delta^m K \left(\max(1,\delta^{-1}\Phi_R+\epsilon)\right)^m
    \Pr(S_{\tau_1} > x\,|\,\tau_1<\infty)\\
    & = \delta^{-1} K \left(\max(\delta,\Phi_R+\delta\epsilon)\right)^m
    \Pr(S_{\tau_1} > x).
  \end{align*}
  It now follows from \eqref{eq:25} that, for sufficiently large $R$,
  we may use the dominated convergence theorem as required to obtain,
  from \eqref{eq:26} and then \eqref{eq:24},
  \begin{align}
    \sum_{m\ge1}\Pr(S_{\tau_m} > x)
    & = (1+o(1))\Pr(S_{\tau_1} > x)\sum_{m\ge1}m\Phi_R^{m-1},
    \qquad\text{as $x\to\infty$,} \nonumber\\
    & \le
    \frac{(1+o(1))}{1-\phg{F}}\overline{F}(x)\sum_{m\ge1}m\Phi_R^{m-1},
    \qquad\text{as $x\to\infty$}.    \label{eq:27}
  \end{align}
  

  Since the random walk~$\{S_n\}_{n\ge0}$ attains its maximum~$M$
  almost surely, for any such sample path, the numbers
  $n=\min\{k\colon{}S_k=M\}$ and $m=\max\{i\colon\tau_i\le{}n\}$ are
  finite.  Then $M=S_{\tau_m}+M^{\tau_m}$ and $M_{\tau_m}<R$.  Hence, for
  $x\ge0$,
  \begin{align}
    \Pr(M>x)
    & \le \sum_{m\ge1}
    \Pr(S_{\tau_m}+M^{\tau_m}>x,\,M^{\tau_m}\in[0,R])
    \nonumber \\
    & = \sum_{m\ge1} \int_0^R \Pr(M^{\tau_m}\in dt,\,S_{\tau_m} > x-t)
    \nonumber \\
    & = \int_0^R \Pr(M\in dt)\sum_{m\ge1}\Pr(S_{\tau_m} > x-t) \label{eq:28}\\
    & \le \frac{(1+o(1))}{1-\phg{F}}\sum_{m\ge1}m\Phi_R^{m-1}
    \int_0^R \Pr(M\in dt)\ov{F}(x-t) \label{eq:29}\\
    & =
    \frac{(1+o(1))}{1-\phg{F}}\overline{F}(x)\sum_{m\ge1}m\Phi_R^{m-1}
    \int_0^R e^{\gamma t}\Pr(M\in dt) \label{eq:30}\\
    & \le
    \frac{(1+o(1))\E e^{\gamma M}}{1-\phg{F}}
    \overline{F}(x)\sum_{m\ge1}m\Phi_R^{m-1}, \label{eq:31}
  \end{align}
  where \eqref{eq:28} follows since, conditional on $\tau_m<\infty$,
  the random variable~$M^{\tau_m}$ is independent of $S_{\tau_m}$ and
  equal in distribution to $M$, the inequality~\eqref{eq:29} follows
  from \eqref{eq:27} (since the integral is taken over the finite
  interval~$[0,R]$), and the equation~\eqref{eq:30} follows from
  \eqref{eq:1}.  Now let $R\to\infty$ in \eqref{eq:31} and use
  \eqref{eq:25} to obtain the required result.
\end{proof}

By combining Lemmas~~\ref{lem:lb} and \ref{lem:ub}, we obtain
Theorem~\ref{thm:sgmax}, where the final equality is given by letting
$\epsilon\to0$.

\begin{theorem}
  \label{thm:sgmax}
  Suppose that $F\in\SS{\gamma}$ for some $\gamma>0$ and that
  $\phg{F}<1$.  Then, given $\epsilon>0$, there exists $a>0$ such that
  \begin{displaymath}
    (1-\epsilon)\frac{\E e^{\gamma M}}{1-\phg{F}}
    \le
    \liminf_{x\to\infty}\frac{1}{\ov{F}(x)}
    \Pr\biggl(M>x,\,\bigcup_{n\ge1}A_n^{a,x-a}\biggr)
    \le
    \lim_{x\to\infty}\frac{\Pr(M>x)}{\ov{F}(x)}
    =
    \frac{\E e^{\gamma M}}{1-\phg{F}}.
  \end{displaymath}
\end{theorem}

\begin{rem}

  It follows in particular from Theorem~\ref{thm:sgmax} that, given
  $\epsilon>0$, there exists $a>0$ such that
  \begin{displaymath}
    \liminf_{x\to\infty}
    \Pr\Biggl(\bigcup_{n\ge1}A_n^{a,x-a} \;\Biggl\lvert\; M > x\Biggr)
    > 1 - \epsilon,
  \end{displaymath}
  that is, for all sufficiently large $x$, conditional on the event
  $\{M>x\}$, at the time~$\tau$ of the first jump of the
  process~$\{S_n\}_{n\ge0}$ above $a$, the value of $S_\tau$ exceeds
  $x-a$ with probability at least $1-\epsilon$.  This is the
  ``principle of a single big jump''.
  We note also that a more compact statement of the conclusion of the
  theorem is that, for any positive function~$h$ satisfying the
  condition~\eqref{eq:5},
  \begin{equation}
    \label{eq:32}
    \lim_{x\to\infty}\frac{1}{\ov{F}(x)}
    \Pr\biggl(M>x,\,\bigcup_{n\ge1}A_n^{h(x),x-h(x)}\biggr)
    =
    \lim_{x\to\infty}\frac{\Pr(M>x)}{\ov{F}(x)}
    =
    \frac{\E e^{\gamma M}}{1-\phg{F}}.
  \end{equation}
  Here again the events $A_n^{h(x),x-h(x)}$, $n\ge1$, are disjoint,
  and we may replace the probability of a union of events by the
  appropriate sum of the individual probabilities.
\end{rem}

\begin{rem}
  While Theorem~\ref{thm:sgmax} does indeed give the asymptotic form
  of $\Pr(M>x)$ as $x\to\infty$, there seems to be no way to determine
  the constant $\E{}e^{\gamma{}M}$,
  which occurs in the statement of the theorem,
  in terms of the distribution function~$F$.  
  Hence the simple bounds given in \eqref{eq:33} may be
  of use in applications.

\end{rem}

We also give a ``spatially local'' result: from
Theorem~\ref{thm:sgmax}, and using again~\eqref{eq:1}, we have
immediately Theorem~\ref{thm:sglocal} below.

\begin{theorem}
  \label{thm:sglocal}
  Suppose that $F\in\SS{\gamma}$ for some $\gamma>0$ and that
  $\phg{F}<1$.  Then, for any fixed $t>0$, and for any positive
  function~$h$ satisfying~\eqref{eq:5},
  \begin{displaymath}
    \lim_{x\to\infty}\frac{1}{\ov{F}(x)}
    \Pr\biggl(M\in(x,x+t],\,\bigcup_{n\ge1}A_n^{h(x),x-h(x)}\biggr)
    =
    \lim_{x\to\infty}\frac{\Pr(M\in(x,x+t])}{\ov{F}(x)}
    =
    \frac{\E e^{\gamma M}}{1-\phg{F}}\left(1 - e^{-\gamma t}\right),
  \end{displaymath}
  with, for any $t_0>0$, uniformity in all $t\in[t_0,\infty]$.
\end{theorem}

\begin{rem}
  The analogous result to Theorem~\ref{thm:sglocal} for the case
  $\gamma=0$ requires the
  slightly stronger condition that $F$ belong to the
  class~$\mathcal{S}^*$ introduced by Kl\"uppelberg (1988)---see
  Asmussen \textit{et al} (2002) and Foss and Zachary (2004).  
  However, for $\gamma>0$ the
  conditions analogous to those defining $\SS{0}$ and
  $\mathcal{S}^*$ match (see also the comments by Rogozin and Sgibnev
  (1999)).
\end{rem}

Analogously to Theorem~\ref{thm:sgmax}, we also have the following
result for the maximum of the random walk on a finite time horizon.

\begin{theorem}
  \label{thm:sgfinite}
  Suppose that $F\in\SS{\gamma}$ for some $\gamma>0$ and that
  $\phg{F}<1$.  Then, for each $N\ge1$, and for any positive
  function~$h$ satisfying~\eqref{eq:5},
  \begin{displaymath}
    \lim_{x\to\infty}\frac{1}{\ov{F}(x)}
    \Pr\biggl(M_N>x,\,\bigcup_{1\le{}n\le{}N}A_n^{h(x),x-h(x)}\biggr)
    =
    \lim_{x\to\infty}\frac{\Pr(M_N>x)}{\ov{F}(x)}
    =
    \sum_{n=1}^N \phg{F}^{n-1}\E e^{\gamma M_{N-n}},
  \end{displaymath}
  where the convergence is uniform over all $N\le\infty$.
\end{theorem}

\begin{proof}
  The proof of Theorem~\ref{thm:sgfinite} is simply a matter of
  checking that the proofs of Lemmas~\ref{lem:lb} and \ref{lem:ub}
  continue to hold (with some small and straightforward adjustments in
  the case of Lemma~\ref{lem:ub}) when restricted to a finite time
  horizon.  The asserted uniformity in $N$ follows from the fact that
  the terms of the form~$\phg{F}^{n-1}$ occurring in \eqref{eq:18} and
  \eqref{eq:23} tend to zero as $n\to\infty$, coupled with the use of a
  simple truncation argument.
\end{proof}
  
\begin{rem}
  \label{rem:bj2}
  Since $\E{}e^{\gamma{}M_n}\to\E{}e^{\gamma{}M}$ as $n\to\infty$, and
  since $\phg{F}<1$, it
  follows from Theorems~\ref{thm:sgmax} and \ref{thm:sgfinite} that,
  given $\epsilon>0$, there exists $N\ge1$ such that
  \begin{displaymath}
    \lim_{x\to\infty}
    \Pr(M_N > x \;\lvert\; M > x)
    > 1 - \epsilon,
  \end{displaymath}
  that is, for all sufficiently large $x$, conditional on the random
  walk ever exceeding $x$, it does so by time~$N$ with probability at
  least $1-\epsilon$.  As remarked in the Introduction, this behaviour
  is very different from that for the case~$F\in\SS{0}$.
\end{rem}

\begin{rem}
  Theorem~\ref{thm:sgfinite} also has a ``spatially local'' version,
  analogous to Theorem~~\ref{thm:sglocal}. 
\end{rem}

The results of Theorem~\ref{thm:sgfinite} may also be extended to
random time horizons (with uniformity in all stopping times
$\sigma\ge1$ a.s.)  We give a result for any a.s.\ finite
stopping time~$\sigma$ such that $S_{\sigma}\le0$ a.s.  
Examples of such stopping times are
\begin{displaymath}
  \sigma_1 = \min\{n\ge 1\colon S_n < 0\} 
  \quad \text {and} \quad
  \sigma_2 = \min\{n>\sigma_1\colon S_n < S_{\sigma_1}\}.
\end{displaymath}
For such a stopping time, let $\chi=-S_{\sigma}\ge0$ and note
that $\Pr(\chi>0)>0$ since $\E\xi_1<0$.

\begin{theorem}\label{thm:sgrandom}
  Suppose that $F\in\SS{\gamma}$ for some $\gamma >0$ and that
  $\phg{F}<1$.  Let $\sigma\ge1$ be an a.s.\ finite stopping time such
  that $S_{\sigma}\le0$ a.s.  Then
  \begin{equation}
    \label{eq:6}
    \lim_{x\to\infty}\frac{\Pr(M_\sigma>x)}{\ov{F}(x)}
    =
    \left(1-\E e^{-\gamma \chi}\right) \, \frac{\E e^{\gamma M}}{1-\phg{F}}.
  \end{equation}
\end{theorem}

\begin{proof}
  Since $\sigma$ is a.s.\ finite and $S_{\sigma}\le0$ a.s., it follows
  from Denisov (2005) that
  \begin{displaymath}
    \Pr(M_{\sigma}>x) \sim \Pr\left(M \in (x,x+\chi')\right),
    \qquad\text{as $x\to\infty$,}
  \end{displaymath}
  where $\chi'$ is a copy of $\chi$ which is independent of $M$.  The
  result~\eqref{eq:6} then follows from Theorem~\ref{thm:sglocal} and
  the dominated convergence theorem.
\end{proof}

For a further extension to any stopping time, one can use the approach
developed by Foss \textit{et al} (2005). Indeed it follows from
Remark~\ref{rem:bj2} that it is sufficient to consider stopping times
which are bounded almost surely from above by a fixed finite number
$N$.  Other extensions are also possible, for example, to ``spatially
local'' versions of these results.

\section*{Acknowledgement}

We are very grateful to Academician A.\ A.\ Borovkov and to Dmitry
Korshunov for helpful and fruitful comments, and to Yuebao Wang for
pointing out a small difficulty in an earlier version of the paper.


\newpage

\begin{thebibliography}{20}



\bibitem{AKKKT} 
  Asmussen, S., Kalashnikov, V., Konstantinides, D., Kl\"uppelberg, C.\ and
  Tsitiashvili, G. (2002).  
  A local limit theorem for random walk maxima with heavy tails.
  \textit{Statist.\ Probab.\ Lett.}, \textbf{56}, 399--404.
  
\bibitem{BD}
  Bertoin, J.\ and Doney, R.\ A. (1996).
  Some asymptotic results for transient random walks.
  {\em Adv.\ Appl.\ Prob.}, {\bf 28}, 207--226.

\bibitem{BOR}
  Borovkov, A.\ A. (1976).
  \emph{Stochastic Processes in Queueing Theory.}
  Springer, New York.

\bibitem{BORBOR} Borovkov, A.\ A.\ and Borovkov, K.\ A. (2004).  
  On probabilities of large deviations for random walks II:  regular
  exponentially decaying distributions.
  \textit{Theory Probab. Appl.}, \textbf{49}, 189--206

\bibitem{CNW}
  Chover, J., Ney, P.\ and Wainger, S. (1973).
  Functions of probability measures. 
  \textit{J. Anal. Math.}, \textbf{26}, 255--302.

\bibitem{DEN} 
  Denisov, D. (2005).
  A Note on the Asymptotics for the Maximum on a Random Time Interval
  of a Random Walk 
  {\em Markov.\  Proc.\ Rel.\ Fields}, \textbf{11}, 165--169.

\bibitem{FPZ}
  Foss, S., Palmowski, Z.\ and Zachary, S. (2005).
  The probability of exceeding a high boundary on a random time
  interval for a heavy-tailed random walk.
  {\em Ann.\ Appl.\ Prob.}, \textbf{15}, 1936--57.

\bibitem{FZ}
  Foss, S.\ and Zachary, S (2003).
  The maximum on a random time interval of a random walk with long-tailed
  increments and negative drift,
  {\em Ann. \ Appl.\ Prob.}, \textbf{13}, 37--53. 

\bibitem{KL88}
  Kl\"uppelberg, C. (1988).
  Subexponential distributions and integrated tails.
  {\em J.\ Appl.\ Prob.}, {\bf 35}, 325--347.
  
\bibitem{PAKES75}
  Pakes, A. (1975).
  On the tails of waiting time distributions.
  {\em J.\ Appl.\ Prob.}, {\bf 7}, 745--789.

\bibitem{PAKES04}
  Pakes, A. (2004).
  Convolution equivalence and infinite divisibility.
  {\em J.\ Appl.\ Prob.}, {\bf 41}, 407--424.

\bibitem{RS}
  Rogozin, B.\ A.\ and Sgibnev, M.\ S. (1999).
  Strongly subexponential distributions and Banach algebras of measures.
  {\em Siberian Math. J.}, {\bf 40}, 963--971.

\bibitem{VER}
  Veraverbeke, N. (1977).
  Asymptotic behavior of Wiener-Hopf factors of a random walk.
  {\em Stoch.\ Proc.\ Appl.}, {\bf 5}, 27--37.
  
\bibitem{ZACH}
  Zachary, S. (2004).
  A note on Veraverbeke's theorem.
  \emph{Queueing Systems}, \textbf{46}, 9--14.

\end{thebibliography}
\end{document}